\documentclass[12pt]{article}
\usepackage{latexsym, amsthm, amsfonts, amsmath, amssymb}

{\theoremstyle{plain}
    \newtheorem{Thm}{\bf Theorem}[section]
    \newtheorem{Prop}[Thm]{\bf Proposition}
    
    \newtheorem{Lma}[Thm]{\bf Lemma}

}
{\theoremstyle{remark}

}
{\theoremstyle{definition}
    \newtheorem{Def}[Thm]{\bf Definition}
}
\numberwithin{equation}{section}

\newcommand{\im}{\operatorname{Im}}

\newcommand{\e}{\operatorname{e}}

\newcommand{\length}{\operatorname{\lambda}}

\newcommand{\fm}{{\mathfrak m}}
\newcommand{\m}{{\mathfrak m}}

\newcommand{\QQ}{{\mathbb{Q}}}

\newcommand{\ringR}{\text{$(R,\fm,k)$ }}

\newcommand{\al}{\alpha}

\newcommand{\brq}{^{[q]}}
\newcommand{\inc}{\subseteq}

\newcounter{hours}\newcounter{minutes}

\newcommand{\excise}[1]{}
\title{\bf When does the F-signature exist?\thanks{  The first author was partially supported by
a grant from the NSA.}}
\author{Ian M. Aberbach\thanks{Department of Mathematics, University of Missouri, Columbia, MO 65211; aberbach@math.missouri.edu} \and Florian Enescu\thanks{Department of Mathematics and Statistics, Georgia State University, Atlanta, 30303 and The Institute of Mathematics of the Romanian Academy, Romania; fenescu@mathstat.gsu.edu}}
\date{}

\begin{document}
\maketitle 

{\bf Abstract.} We show that the $F$-signature of an $F$-finite local ring $R$ of characteristic $p >0$ exists when $R$ is either the localization of an $\mathbf{N}$-graded ring at its irrelevant ideal or $\mathbf{Q}$-Gorenstein on its punctured spectrum. This extends results by Huneke, Leuschke, Yao and Singh and proves the existence of the $F$-signature in the cases where weak $F$-regularity is known to be equivalent to strong $F$-regularity.

{\bf R\'esum\'e.} Nous prouvons dans cet article l'existence de la F-signature d'un
anneau local F-fini R, de caract\'eristique positive p, quand R est la localisation \`a
l'unique id\'eal homog\`ene maximal d'un anneau $\mathbf{N}$-gradu\'e ou quand R est $\mathbf{Q}$-Gorenstein sur son spectre \'epoint\'e. Ceci g\'en\'eralise les r\'esultats de Huneke, Leuschke, Yao et Singh et
prouve l'existence de la F-signature dans les cas o\`u faible et forte F-r\'egularit\'e sont \'equivalentes.

\section{A sufficient condition for the existence of the F-signature}

Let $\ringR$ be a reduced, local $F$-finite
ring of positive characteristic $p>0$
and Krull dimension $d$.  Let $$R^{1/q} = R^{a_q} \oplus M_q$$ 
be a
direct sum decomposition of $R^{1/q}$ such that $M_q$ has no free
direct summands. If $R$ is complete, such a decomposition is unique up
to isomorphism. Recent research has focused on the asymptotic growth rate of
the numbers $a_q$ as $q \to \infty$.  In particular, the $F$-signature
(defined below) is studied in~\cite{HL} and~\cite{AL}, and more generally
the Frobenius splitting ratio is studied in~\cite{AE}.

For a local ring $\ringR$,
we set $\alpha(R) =\log _p[k_R:k_R^p]$. It is easy to see that, for
an $\fm$-primary ideal $I$ of $R$, $\length (R^{1/q}/IR^{1/q})=
\length(R/I^{[q]})/ q^{\al(R)}$, where $\length(-)$ represents the length function over $R$. 

We would like to first define the notion of $F$-signature as it appears in~\cite{AL} and~\cite{HL}.

\begin{Def}
The {\sl $F$-signature} of $R$  
is $s(R) = \lim_{q \to \infty} \dfrac {a_q}{q^{d+\al(R)}}$, if it exists.
\end{Def}

The following result, due to Aberbach and Leuschke~\cite{AL}, holds:

\begin{Thm}
Let $\ringR$ be a reduced Noetherian ring of positive characteristic $p$. Then
$\liminf_{q \to \infty} a_q/q^{d+\al(R)} >0$ if and only if $\limsup_{q \to \infty} a_q/q^{d+\al(R)} >0$ if and only if $R$ is strongly $F$-regular.
\end{Thm}

The question of whether or not, in a strongly $F$-regular ring, $s(R)$ exists,
is open.  We show in this paper that its existence is closely connected
to the question of whether or not weak and strong $F$-regularity are
equivalent.

Smith and Van den Bergh (\cite{SV}) have shown that the $F$-signature of $R$ exists when  $R$ has finite Frobenius representation type (FFRT) type, that is, if only finitely many isomorphism classes of indecomposable maximal Cohen-Macaulay modules occur as direct summands of $R^{1/q}$ for any $q=p^e$. Yao has proven that, under mild conditions, tight closure commutes with localization in a ring of FFRT type, \cite{Y1}. Moreover, Huneke and Leuschke proved that if $R$ is also Gorenstein, then the $F$-signature exists, \cite{HL}. Yao has recently extended this result to rings that are Gorenstein on their punctured spectrum, \cite{Y2}. Singh has also shown that the $F$-signature exists for monomial rings, \cite{Si}. 

Let $(R,\fm)$ be an approximately Gorenstein ring.  This means that
$R$ has a sequence of $\fm$-primary irreducible ideals $\{I_t\}_t$ cofinal
with the powers of $\fm$.  By taking a subsequence, we may assume that
$I_t \supset I_{t+1}$.  For each $t$, let $u_t$ be an element of $R$ which
represents a socle element modulo $I_t$.  Then there is, for each $t$, a
homomorphism $R/I_t \hookrightarrow R/I_{t+1}$ such that $u_t +I_t\mapsto u_{t+1} + I_{t+1}$.  The direct limit of the system will be  the injective hull $E =E_R(R/\fm)$
and each $u_t$ will map to the socle element of $E$, which we will denote by $u$.  Hochster has shown that every excellent, reduced local ring is
approximately Gorenstein (\cite{Ho}).  

Aberbach and Leuschke have shown that, for every $q$, there exists $t_0 (q)$, such that 

$$ a_q/(q^{d+ \alpha(R)}) = \length(R/(I_t\brq  : u_t^q))/ q^d,$$ 
for all $t \geq t_0(q)$ (see~\cite{AL}, p.~55).

The situation when $t_0(q)$ can be chosen independently of $q$ is of special interest.

\begin{Def}
\label{A}
We say that $R$ satisfies {\sl Condition $(A)$}, if there exist 
a sequence of irreducible $\m$-primary ideals $\{I_t\}$
and a $t_0$ such that, for all $t \geq t_0$ and all $q$
$$(I_t\brq  : u_t^q) = (I_{t_0}\brq  : u_{t_0}^q).$$
\end{Def}

\begin{Prop}
Let $\ringR$ be a local reduced $F$-finite ring. If $R$ satisfies Condition $A$, then the $F$-signature exists.
\end{Prop}

\begin{proof}
We know that $R$ is approximately Gorenstein and hence we will use the notation fixed in the paragraph above. 

As explained above, Condition $A$ implies that there exists $t_0$, independent of $q$, such that 
$$ a_q/(q^{d+ \alpha(R)}) = \length(R/(I_{t_0}\brq  : u_{t_0}^q))/ q^d,$$ for all $q$.

But $\length(R/(I_{t_0}\brq  : u_{t_0}^q))= 
\length(R/I_{t_0}\brq) -\length(R/(I_{t_0} + u_{t_0}R)\brq)$.
Dividing by $q^d$ and taking the limit as $q \to \infty$ yields
  $s(R) = \e_{HK}(I_{t_0}, R) - \e_{HK}(I_{t_0}+u_{t_0}R, R)$. 
\end{proof}

Now we would like to concentrate on another condition, Condition $(B)$, 
that appeared first in the work of Yao. First we need to introduce
 some notation.

Assume that $E$ is the injective hull of the residue field $k$. By $R^{(e)}$ we denote the $R$-bialgebra whose underlying abelian group equals $R$ and the left and right $R$-multiplication is given by $ a \cdot r * b = arb^q$, for $a,b \in R, r \in R^{(e)}$.

Let $k = Ru \to E$ be the natural inclusion and consider the natural induced map $\phi _e: R^{(e)} \otimes _R E \to R^{(e)} \otimes_R (E/k)$. Then $a_q / q ^{\alpha(R)}  = \length (\ker (\phi_e))$ (by Aberbach-Enescu, Corollary 2.8 in~\cite{AE}, see also Yao's work~\cite{Y2}).

One can in fact see that $$\length (\ker (\phi_e)) = \length (R/( c \in R: c \otimes u = 0  \textrm{ in }  R^{(e)} \otimes_R E )) = \length (R/\cup_t (I_t \brq : u_t^q)).$$

\begin{Def}
\label{B}
We say that $R$ satisfies {\sl Condition $(B)$} if there exists a finite length submodule $E' \subset E$ such that, if $\psi _e : R^{(e)} \otimes_R E' \to R^{(e)} \otimes_R E'/k$, then $\length (\ker (\phi_e)) = \length (\ker (\psi_e))$, for all $e$.
\end{Def}

Yao~\cite{Y2} has shown that Condition $(B)$ implies that the $F$-signature of $R$ exists.

\begin{Prop}
Let $\ringR$ be a local reduced $F$-finite ring. Then Conditions $(A)$ and $(B)$ are equivalent.
\end{Prop}

\begin{proof}
Assume that  Condition $(A)$ holds. Then one can take $E' = R/I_{t_0}$ and then Condition $(B)$ follows.

If Condition $(B)$ holds, then take $t_0$ large enough such that 
$E' \subset \im (R/I_{t_0} \to E)$.

As noted above, one can compute the length of the kernel of $\psi_e$ as the colength of $\{ c \in R: c \otimes u = 0 \textrm{ in } R^{(e)} \otimes_R E' \}$. Since $R/I_{t_0}$ injects into $E$ we see that  $\{ c \in R: c \otimes u = 0  \textrm{ in }  R^{(e)} \otimes_R E' \}$ is a subset of $\{ c \in R: c \otimes u = 0  \textrm{ in } R^{(e)} \otimes_R R/I_{t_0} \} = (I_{t_0}\brq : u_{t_0}^q ) $.

Since $(I_{t_0}\brq : u_{t_0}^q ) \subset (I_{t}\brq : u_{t}^q ) $ for all $t \geq t_0$, we see that Condition $(B)$ implies that $(I_{t_0}\brq : u_{t_0}^q ) = (I_{t}\brq : u_{t}^q ) $ for all $t \geq t_0$, which is Condition $(A)$.

\end{proof}

\section{$\mathbf{N}$-Graded Rings}

Let $(R, \fm)$ be a Noetherian $\mathbf{N}$-graded ring $R= \oplus_{n \geq 0} R_n$, where $R_0=k$ is an $F$-finite field of characteristic $p >0$.

For any graded $R$-module $M$ one can define a natural grading on $R^{(e)} \otimes M$: the degree of any tensor monomial $r \otimes m$ equals 
$\deg(r) + q\deg(m).$

In what follows we will need the following important Lemma by Lyubeznik and Smith (\cite{LS}, Theorem 3.2):

\begin{Lma}
\label{LS}
Let $R$ be an $\mathbf{N}$-graded ring and $M, N$ two graded $R$-modules. Then there exists an integer $t$ depending only on $R$ such that
whenever $$M \to N$$ is a degree preserving map which is bijective in degrees greater than $s$, then the induced map $$ R^{(e)} \otimes M \to R^{(e)} \otimes N$$ is bijective in degrees greater than $p^e(s+t)$.

\end{Lma}

Let $E$ be the injective hull of $R_{\fm}$. In fact, $E$ is also the injective hull of $R/\fm$ over $R$ and as a result is naturally graded with socle in
degree $0$. We can write $E = \oplus_{n \leq 0} E_n$.

Let $t$ be as in the Lemma ~\ref{LS}, and let $s \leq  -t-1 $. Obviously the map $E' = \oplus_{s \leq n \leq 0} E_n \to E = \oplus_{n \leq 0} E_n$ is bijective in degrees greater than $s$. So by  Lemma~\ref{LS}, the map $R^{(e)} \otimes E' \to R^{(e)} \otimes E$ is bijective in degrees greater than $p^e(s+t) \leq -p^e$.

\begin{Thm}
Let $R$ be an $\mathbf{N}$-graded reduced ring over an $F$-finite field $k$ of positive characteristic. Then Condition $(B)$ is satisfied by $R$ and hence the $F$-signature of $R$ exists.
\end{Thm}

\begin{proof}
Let $E$ be the injective hull of $k= R/\fm$ over $R_\fm$. As above, $E = \oplus_{n \leq 0} E_n$, where $0$ is the degree of the socle generator $u$ of $E$.

Using the notation introduced above, we will let $s = -t-1$ and $E'  = \oplus_{s \leq n \leq n_0} E_n \to E$. So, $R^{(e)} \otimes E' \to R^{(e)} \otimes E$ is bijective in degrees greater than $-p^e$. In particular it is bijective in degrees greater or equal to $0$.

We have the following exact sequences:

$$ 0 \to k = Ru \to E \to E/k \to 0$$ and

$$ 0 \to k = Ru \to E' \to E'/k \to 0.$$

After tensoring with $R^{(e)}$, we get the exact sequences

$$ R^{(e)} \otimes k = R^{(e)} \otimes Ru \to R^{(e)} \otimes E \stackrel{\phi_e}{\to} R^{(e)} \otimes E/k \to 0$$ and

$$ R^{(e)} \otimes k = R^{(e)} \otimes Ru \to R^{(e)} \otimes E' \stackrel{\psi_e}{\to} R^{(e)} \otimes E'/k \to 0.$$

One can easily see that $\ker(\phi_e)$ and $\ker(\psi_e)$ are the
submodules generated by $1 \otimes u$  in $R^{(e)} \otimes E$ and $R^{(e)} \otimes E'$, respectively.

The degree of $1 \otimes u$ is $q\cdot 0 =0$ and we have noted that the natural map $R^{(e)} \otimes E' \to R^{(e)} \otimes E$ is bijective in degrees greater than $-p^e$. This shows that $\ker(\phi_e) \simeq \ker(\psi_e)$ and hence Condition $(B)$ is satisfied.

\end{proof}

\section{$\mathbb{Q}$-Gorenstein Rings}

We turn now to showing that Condition (A) holds in strongly $F$-regular 
local rings which are $\QQ$-Gorenstein on the punctured spectrum.
Let $\ringR$ be such a ring of dimension $d$, and assume that $R$ has
a canonical module (e.g. $R$ is complete).  In this case $R$ has an unmixed
ideal of height 1, say $J \subseteq R$, which is a canonical ideal.
We may pick an element $a \in J$ which generates $J$ at all minimal primes
of $J$, and then an element $x_2 \in \m$ which is a parameter on $R/J$
such that $x_2 J \inc aR$.  It is easy to see that then $x^n J^{(n)} \inc
a^n R$ for all $n \ge 1$  (where  $J^{(n)}$ is the height one component of
$J^n$).  The condition that $R$ is $\QQ$-Gorenstein on the punctured spectrum
implies that there is an integer $h$ and two sequences of elements
$x_3,\ldots, x_d \in \m$ and $a_3,\ldots, a_d \in J^{(h)}$ such that
$x_i J^{(h)} \inc a_iR$ for $3 \le i \le d$, and $x_2, \ldots, x_d$ is a
s.o.p.~on $R/J$.  We may then pick $x_1 \in J$ such that $x_1,\ldots, x_d$
is an s.o.p.~for $R$.  See \cite{Ab}, section 2.2 for more detail.  Then
by \cite{Ab}, Lemma 2.2.3 we have that for any $N \ge 0$ and any $n \ge 0$,
\begin{equation}\label{colon}
(J^{(nh)}, x_2^N, \ldots, \widehat{x_i^N},\ldots, x_d^N):x_i^\infty
= (J^{(nh)}, x_2^N, \ldots, \widehat{x_i^N},\ldots, x_d^N):x_i^n.
\end{equation}

\begin{Thm} Let $\ringR$ be an $F$-finite strongly $F$-regular ring which is
$\QQ$-Gorenstein on the punctured spectrum.  Then $R$ satisfies Condition
{\rm (A)}.  In particular the $F$-signature of $R$ exists.
\end{Thm}

\begin{proof}
If $R$ is not complete, we observe that, since $R$ is excellent, 
$\widehat R$ is strongly $F$-regular and
$\QQ$-Gorenstein on the punctured spectrum.   If $\{I_t\}$ is a sequence of ideal in $\widehat R$ showing condition
(A) in $\widehat R$, then $\{ I_t \cap R\}$ does so for $R$.  Thus we will assume
that $R$ is complete.

Let $J$, $h$, and $x_1,\ldots, x_d$ be as discussed above.  Let
$I_t = (x_1^{t-1}J, x_2^t, \ldots, x_d^t)$.  
Since $x_1^n J \cong J$ as $R$-modules, the quotient $R/x_1^nJ$ is Gorenstein.
The hypothesis that $x_2,\ldots, x_d$ are parameters on $R/J$ and $R/x_1R$
(hence on $R/x_1^nJ$) then shows that $I_t$ is irreducible (see \cite{BH},
Proposition 3.3.18).
The sequence $\{I_t\}$ is then
a sequence of $\m$-primary irreducible ideals cofinal with the powers
of $\m$.  If $u_1$ represents the socle element of $I_1$, then 
we may take $u_t = (x_1\cdots x_d)^{t-1}u_1$ to represent the socle 
element of $I_t$.  We will show that $t_0$ may be taken to be $3$.

Suppose that $c \in I_t\brq: u_t^q$ for some $q$.  We will show that
$c \in I_3\brq: u_3^q$.  Raising to the $q'$th power we have
$c^{q'} u_t^{qq'} = c^{q'} \left( (x_1\cdots x_d)^{t-1}u_1\right)^{qq'}
\in I_t^{[qq']} = ( x_1^{t-1} J, x_2^t,\ldots, x_d^t)^{[qq']}$.
Hence $c^{q'}\left((x_2\cdots x_d)^{t-1}u_1\right)^{qq'}
\in (x_2^t,\ldots, x_d^t)^{[qq']} : x_1^{(t-1)qq'} + 
(J, x_2^t,\ldots, x_d^t)^{[qq']} = (J, x_2^t,\ldots, x_d^t)^{[qq']}$.

Write $qq' = n_{q'} h + r_{q'}$ with $0 \le r_{q'} < h$.  Repeated application
of equation \ref{colon} (using $1$ rather than $h$ for $x_2$) gives
\begin{equation}
c^{q'}( (x_2\cdots x_d)u_1)^{qq'} \in (J^{(n_{q'}h)}, x_2^{2qq'},
\ldots, x_d^{2qq'}).
\end{equation}
Let $d \in J^{(h)} \inc J^{(r_{q'})}$.
Multiplying by $x_2^{qq'}$ and using that $x_2^{qq'}J^{(qq')} \inc a^{qq'}R
\inc J^{[qq']}$ we have $ dc^{q'}\left((x_2 \cdots x_d)^2u_1\right)^{qq'}
\in (J, x_2^3,\ldots, x_d^3)^{[qq']}$.   Multiplying by $x_1^{2qq'}$ shows
that $dc^{q'}u_3^{qq'} = d(cu_3^q)^{q'} \in (I_3\brq)^{[q']}$.  Thus
$cu_3^q \in (I_3\brq)^* = I_3\brq$, as desired.
\end{proof}

\end{document}